\documentclass[onecolumn,showpacs,showkeys,preprintnumbers]{revtex4}

\usepackage{amsmath}
\usepackage{amssymb}
\usepackage{graphicx}
\usepackage{dcolumn}
\usepackage{bm}
\usepackage{mathrsfs}
\usepackage{enumerate}
\usepackage{hyperref}

\newtheorem{theorem}{Theorem}[section]
\newtheorem{lemma}[theorem]{Lemma}

\DeclareMathAlphabet{\mathpzc}{OT1}{pzc}{m}{it}

\begin{document}
\title{Controllability of Quantum Systems on the Lie Group $SU(1,1)$}
\thanks{This research was supported in part by the National Natural
Science Foundation of China under Grant Number 60674039 and
60433050. Tzyh-Jong~Tarn would also like to acknowledge partial
support from the U.S. Army Research Office under Grant
W911NF-04-1-0386.}

\author{Jian-Wu Wu$^{1}$} \email{wujw03@mails.tsinghua.edu.cn}
\author{Chun-Wen Li$^{1}$}
\author{Jing Zhang$^{1}$}
\author {Tzyh-Jong Tarn$^{2}$}

\affiliation{$^{1}$Department of Automation, Tsinghua University,
Beijing, 100084, P. R. China\\
$^{2}$Department of Systems Science and Mathematics, Washington
University, St. Louis, MO 63130, USA}

\date{\today}

\begin{abstract}
This paper examines the controllability for quantum control
systems with $SU(1,1)$ dynamical symmetry, namely, the ability to
use some electromagnetic field to redirect the quantum system
toward a desired evolution. The problem is formalized as the
control of a right invariant bilinear system evolving on the Lie
group $SU(1,1)$ of two dimensional special pseudo-unitary
matrices. It is proved that the elliptic condition of the total
Hamiltonian is both sufficient and necessary for the
controllability. Conditions are also given for small time local
controllability and strong controllability. The results obtained
are also valid for the control systems on the Lie groups $SO(2,1)$
and $SL(2,\mathbb{R})$.
\end{abstract}



\maketitle

\section{Introduction}\label{sec1}
Controllability is a fundamental problem in the control theory
with respect to both classical
\cite{Sussmann1,Brockket1,Jurdjevic1,Sussmann2} and
quantum~\cite{Huang1,Ramakrishna1,Tarn1,Alessandro1,Schirmer1,Lan1,wrb1}
mechanical system. In the past decades, sufficient conditions
\cite{Jurdjevic1,Huang1,Ramakrishna1,Lan1,wrb1} have been
established via algebraic methods for systems evolving on
manifolds or Lie groups. However, most of these conditions are not
necessary, especially for the systems on noncompact Lie
groups\cite{{Jurdjevic1}}.

The main purpose of this article is to establish a sufficient and
necessary condition that examines the controllability of the
quantum systems whose propagators evolve on the noncompact Lie
group $SU(1,1)$, which describes the dynamical symmetry of many
important physical possesses, e.g., the downconversion
process~\cite{Puri1,Gerry2}, the Bose-Einstein
condensation~\cite{Gortel1}, the spin wave transition in
solid-state physics~\cite{Bose1}, the evolution in free
space~\cite{Agarwal1}.

The problem is investigated by considering the following right
invariant bilinear system on the Lie group $SU(1,1)$%
\noindent\begin{equation}\label{eq1.1}
\dot{X}(t)=\left[A+\sum\limits_{l=1}^{r}u_l(t)B_l\right]X(t),
X(0)=I_2,
\end{equation}

\noindent where $u_l(t)$ belong to some admissible control set
$\mathcal{U}$, which consists of functions defined on
$\mathbb{R}^+=[0,\infty)$. The drift term $A$ and the control
terms $B_1$, $B_2$, $\cdots$, $B_r$ are elements of the Lie
algebra $su(1,1)$, where $B_1$, $B_2$, $\cdots$, $B_r$ are assumed
to be linearly independent with respect to real coefficients. The
state, $X(t)$, is a two-dimensional complex pseudo-unitary
matrix in the form of%
\noindent\begin{equation}\label{eq1.2}
 \left(\begin{array}{cc}
  a & b \\
  \bar{b} & \bar{a}\\
\end{array}\right),~~~~~|a|^2-|b|^2=1,
\end{equation}

\noindent where $\bar{a}$ represents the complex conjugate of $a$.
Since $SU(1,1)$ is homomorphic to $SO(2,1)$ and isomorphic to
$SL(2,\mathbb{R})$ respectively, the results obtained in this
paper are still valid for the systems on these two Lie groups.

For a driftless system varying on the noncompact Lie group, it was
shown in~\cite{Jurdjevic1} that the system is controllable when
there exists a constant control such that the state trajectory is
periodic. Applied to the quantum system evolving on $SU(1,1)$, it
can be concluded that the system is controllable if the total
Hamiltonian (including the internal Hamiltonian and the
interaction Hamiltonian) of the system can be adjusted to be
elliptic. In \cite{wjw1}, this sufficient condition was extended
to bounded controls, algorithms were given accordingly to design
control laws to achieve desired evolutions. In this paper, this
condition is proven to be necessary for the single input case,
which can be directly used to judge whether one can find a
``magnetic field'' to induce a desired transition between two
arbitrary $SU(1,1)$ coherent states, which are of particular
importance in quantum optics~\cite{Gerry1,Walls1}.

The paper is organized as follows. Section~\ref{sec2} presents
preliminaries to be used in the rest of this paper.
Subsection~\ref{sec2.1} describes the systems to be considered in
mathematical terms of right invariant bilinear control systems
that evolve on the noncompact Lie group $SU(1,1)$.
Subsection~\ref{sec2.2} introduces necessary definitions for the
system controllability. Section~\ref{sec3} contains the main
results on the system controllability. In Subsection~\ref{sec3.1},
we present some properties of Lie algebra $su(1,1)$ that will be
useful for studying the system controllability. In
Subsection~\ref{sec3.2}, a sufficient and necessary condition that
examines the controllability is established for the single input
case, showing that the controllability of such quantum systems can
be completely determined by finding a constant control that
adjusts the total Hamiltonian of the undergoing system to be
elliptic, or not. Properties of the strong controllability and
small time local controllability are discussed in the subsequence
as well. Controllability properties for the multi-input case is
considered in Subsection~\ref{sec3.3}. In Section~\ref{sec4}, we
discuss the relationship between systems evolving on $SO(2,1)$,
$SL(2,\mathbb{R})$ and $SU(1,1)$ and show that the result obtained
are still valid for the system evolving on these noncompact Lie
groups. Interpretations are then provided for the criteria
obtained based on the topology of $SO(2,1)$. Illustrative examples
are elaborated in Section~\ref{sec6}. Finally, conclusions are
drawn in section~\ref{sec6}.

\section{Preliminaries}\label{sec2}

In this section, we present preliminaries which will be used in
this paper.

\subsection{Quantum Control Systems on $SU(1,1)$}\label{sec2.1}

The time evolution of a controlled quantum system is determined
through the Schr\"{o}dinger equation%
\begin{equation}\label{eq2.1.1}
i\hbar\frac{d}{dt}\psi(t)=\left[H_0+\sum_{l=1}^{r}u_l(t)H_l\right]\psi(t),~\psi(0)=\psi_0,
\end{equation}

\noindent where the wave function $\psi(t)$ describes the state of
the system in an appropriate Hilbert space $\mathcal{H}$. The
Hermitian operators $H_0$ and $H_l~(l=1,2,\cdots,r)$ are referred
to as the internal and interaction Hamiltonians respectively. The
scalars $u_l(t)~(l=1,2,\cdots,r)$ represent some adjustable
classical fields coupled to the system, which are used to control
the evolution of the system.

In this paper, we study the class of quantum systems evolving on
the noncompact Lie group $SU(1,1)$, whose internal and interaction
Hamiltonians can be expressed as linear combinations of the
operators $K_x$, $K_y$ and $K_z$, which satisfy the following
commutation
relations%
\noindent\begin{equation}\label{eq2.1.2}
[K_x,K_y]=-iK_z,~[K_y,K_z]=iK_x,~[K_z,K_x]=iK_y,
\end{equation}

\noindent i.e., closed as an $su(1,1)$ Lie algebra. According to
the group representation theory~\cite{Vilenkin1}, $H_0$ and
$H_l~(l=1,2,\cdots,r)$ are all operators on an infinite
dimensional Hilbert space $\mathcal{H}$ because $su(1,1)$ is
noncompact (see Example~1).

Let $U(t)$ be the evolution operator~(or propagator) that
transforms the system state from the initial $\psi(0)$ to
$\psi(t)$, i.e., $\psi(t)=U(t)\psi(0)$. Then, from
(\ref{eq2.1.1}), by setting $\hbar=1$ and
$\bar{H}_l=-iH_l~(l=0,1,\cdots,r)$, we can obtain that%
\noindent\begin{equation}\label{eq2.1.3}
\dot{U}(t)=\left[\bar{H}_0+\sum_{l=1}^{r}u_l(t)\bar{H}_l\right]U(t),~U(0)=I,
\end{equation}

\noindent where $I$ is the identity operator on $\mathcal{H}$. The
evolution operator $U(t)$ can be treated as an infinite
dimensional matrix since it acts on the infinite dimensional
states space. It is inconvenient to study the controllability
properties of such infinite-dimensional systems directly.
Nevertheless, since all faithful representations are algebraically
isomorphic on which the system controllability property does not
rely, one can always focus the study on the equivalent system
(\ref{eq1.1}) evolving on the Lie group $SU(1,1)$ of
pseudo-unitary matrices, where $A$ and $B_l$ can be written down
as linear combinations of
\noindent\begin{equation}\label{eq2.1.4}%
\begin{array}{l}
 \bar{K}_x=\frac{1}{2}\sigma_y=\frac{1}{2}\left(\begin{array}{cc}
  0 & -i \\
  i & 0\\
\end{array}\right),~~
\bar{K}_y=-\frac{1}{2}\sigma_x=\frac{1}{2}\left(\begin{array}{cc}
  0 & -1 \\
  -1 & 0\\
\end{array}\right),~~
 \bar{K}_z=-\frac{i}{2}\sigma_z=\frac{1}{2}\left(\begin{array}{cc}
  -i & 0 \\
  0 & i\\
\end{array}\right),
\end{array}
\end{equation}

\noindent where $\sigma_{x,y,z}$ are Pauli matrices. The matrices
$\bar{K}_x$, $\bar{K}_y$ and $\bar{K}_z$ are non-unitary
representation of the operators $K_x$, $K_y$ and $K_z$, and one
can verify that $\bar{K}_x$, $\bar{K}_y$ and $\bar{K}_z$ satisfy%
\noindent\begin{equation}\label{eq2.1.5}
[\bar{K}_x,\bar{K}_y]=-\bar{K}_z,~[\bar{K}_y,\bar{K}_z]=\bar{K}_x,~[\bar{K}_z,\bar{K}_x]=\bar{K}_y.
\end{equation}

\noindent $\bar{K}_x$, $\bar{K}_y$ and $\bar{K}_z$ form a
normalized orthonormal basis of the Lie algebra $su(1,1)$ with
respect to the inner product
$\left<\cdot,\cdot\right>$ defined by%
\noindent\begin{equation}\label{eq2.1.6}
\left<M,~N\right>=2~\mbox{Tr}(MN^\dag),
\end{equation}

\noindent where $N^\dag$ denotes the Hermitian conjugation of $N$.
As a result, any given element $M$ in $su(1,1)$ can be expressed
in the
following way%
\noindent\begin{equation}\label{eq2.1.7}
M=\left<M,\bar{K}_x\right>\bar{K}_x+\left<M,\bar{K}_y\right>\bar{K}_y+\left<M,\bar{K}_z\right>\bar{K}_z.
\end{equation}

\subsection{Controllability and the Reachable Sets}\label{sec2.2}

To define the controllability of system (\ref{eq1.1}), the
following reachable sets started from the identity $I_2$ are
useful,

\noindent\begin{enumerate}%
    \item $R(t)=\{X_f|\exists~u, s.t., X(0)=I_2, X(t)=X_f\}$, i.e., the set of all the possible values for state
    $X(t)$ at time $t$.
    \item $R(\cup{T})=\bigcup_{0\leq{t}\leq{T}}R(t)$
    ($R(\cup{\infty})=\bigcup_{0\leq{t}<\infty}R(t)$), i.e.,
    the set of all the possible states
    within time $T$ ($\infty$).
    \item $R(\cap{T})=\bigcap_{0<{t}\leq{T}}R(t)$
    ($R(\cap{\infty})=\bigcap_{0<{t}<\infty}R(t)$), i.e.,
    the set of all values for state $X(t)$ that can be achieved
    at any time within $T$ ($\infty$).
\end{enumerate}

With the reachable sets defined above, controllability of the
system (\ref{eq1.1}) can be defined as follows.

\noindent{\bf{Definition~1:}} System (\ref{eq1.1}) is said to be
controllable on $SU(1,1)$ if $R(\cup{\infty})=SU(1,1)$, strongly
controllable on $SU(1,1)$ if $R(\cap{\infty})=SU(1,1)$, and small
time local controllable on $SU(1,1)$ if $I_2$ is an interior point
of $R(t)$ for any $t>0$.

Apparently, system (\ref{eq1.1}) is both controllable and small
time local controllable if it is strongly controllable. The right
invariant property indicates that the controllability properties
of system (\ref{eq1.1}) is independent of the system initial
condition.

Before discussing the controllability of system (\ref{eq1.1}), we
introduce the following three Lie algebras.

\noindent\begin{enumerate}%
\item $\mathcal{L}$ is the Lie algebra
generated by $\{A,B_1,B_2,\cdots,B_r\}$, and $e^\mathcal{L}$ is
the connected Lie subgroup of $SU(1,1)$ exponentiated by
$\mathcal{L}$.

\item $\mathcal{L}_0$ is the maximal ideal in $\mathcal{L}$
generated by $\{B_1,B_2,\cdots,B_r\}$, and $e^{\mathcal{L}_0}$ is
the connected Lie subgroup of $SU(1,1)$ exponentiated by
$\mathcal{L}_0$.

\item $\mathcal{B}$ is the algebra generated by
$\{B_1,B_2,\cdots,B_r\}$, and $e^\mathcal{B}$ is the connected Lie
subgroup of $SU(1,1)$ exponentiated by $\mathcal{B}$.
\end{enumerate}

Clearly, $R(\cup{\infty}){\subseteq}e^\mathcal{L}$, which implies
that $e^\mathcal{L}$ must equal $SU(1,1)$ if system (\ref{eq1.1})
is controllable. $\mathcal{L}_0$ has co-dimension 0 or 1 in
$\mathcal{L}$ depending on whether $A$ is an element of
$\mathcal{L}_0$ or not.

\section{Controllability}\label{sec3}

In this section, the results on the controllability of the systems
with respect to both single-input and multi-input cases will be
presented. For that purpose, the following properties of Lie
algebra $su(1,1)$ will be very useful.

\subsection{Properties of Lie Algebra $su(1,1)$}\label{sec3.1}

{\bf{Definition~2:}} An $su(1,1)$ element $M$ (as well as its
exponential $\exp(tM)$) is called elliptic (hyperbolic, parabolic)
if $\left<M,M^\dag\right>$ is negative (positive, zero).

\noindent\begin{lemma}\label{lemma3.1.1}%
The commutator $[M,N]$ is elliptic (parabolic, or hyperbolic) if
and only if
$\left<M,N^{\dag}\right>^2-\left<N,N^{\dag}\right>\left<M,M^{\dag}\right><0$
($=0$, or $>0$).
\end{lemma}

{\it{Proof:~}} Since $\bar{K}_x$, $\bar{K}_y$ and
$\bar{K}_z$ span the Lie algebra $su(1,1)$, we can write%
\noindent\begin{equation}\label{eq3.1.1}%
M=m_1\bar{K}_x+m_2\bar{K}_y+m_3\bar{K}_z,
\end{equation}

and%
\noindent\begin{equation}\label{eq3.1.2}%
N=n_1\bar{K}_x+n_2\bar{K}_y+n_3\bar{K}_z,
\end{equation}

\noindent where the coefficients $m_l$ and $n_l$ are real numbers.
Making use of the commutation relations given in
(\ref{eq2.1.2}), we have%
\noindent\begin{equation}\label{eq3.1.3}
  [M,N]=(m_2n_3-m_3n_2)\bar{K}_x+(m_3n_1-m_1n_3)\bar{K}_y-(m_1n_2-m_2n_1)\bar{K}_z.
\end{equation}

\noindent A simple computation yields that%
\noindent\begin{equation}\label{eq3.1.4}%
\begin{array}{l}
   \left<M,M^{\dag}\right>=m_1^2+m_2^2-m_3^2,~~
   \left<N,N^{\dag}\right>=n_1^2+n_2^2-n_3^2,~~
   \left<M,N^{\dag}\right>=m_1n_1+m_2n_2-m_3n_3,
\end{array}
\end{equation}

\noindent and%
\noindent\begin{equation}\label{eq3.1.5}%
\left<[M,N],[M,N]^\dag\right>=(m_2n_3-m_3n_2)^2+(m_3n_1-m_1n_3)^2-(m_1n_2-m_2n_1)^2.
\end{equation}

\noindent Comparison of (\ref{eq3.1.4}) and (\ref{eq3.1.5}) gives%
\noindent\begin{equation}\label{eq3.1.6}
\left<[M,N],[M,N]^{\dag}\right>=\left<M,N^{\dag}\right>^2-\left<N,N^{\dag}\right>\left<M,M^{\dag}\right>.
\end{equation}

\noindent The statement of the Lemma follows immediately from the
above equation.~~~~~~~~~~~~~~~~~~~~~~~~~~~~~~~~~~~~~$\square$

\noindent\begin{lemma}\label{lemma3.1.2}%
Given any two linearly independent elements $M$ and $N$ in
$su(1,1)$, $M$, $N$ and $[M,N]$ are linearly independent if and
only if $[M,N]$ is not parabolic.
\end{lemma}

{\it{Proof:~}} From (\ref{eq3.1.1})-(\ref{eq3.1.3}), it can be
concluded that $M$, $N$ and $[M,N]$ are linearly
independent if and only if%
\noindent\begin{equation}\label{eq3.1.7}
\left|%
\begin{array}{ccc}
  m_1 & n_1 & m_2n_3-m_3n_2 \\
  m_2 & n_2 & m_3n_1-m_1n_3 \\
  m_3 & n_3 & -(m_1n_2-m_2n_1)\\
\end{array}%
\right|\neq0,
\end{equation}

\noindent or equivalently%
\noindent\begin{equation}\label{eq3.1.8}
(m_2n_3-m_3n_2)^2+(m_3n_1-m_1n_3)^2-(m_1n_2-m_2n_1)^2\neq0,
\end{equation}

\noindent i.e., $\left<[M,N],[M,N]^\dag\right>\neq0$. It
immediately follows from Eq.(\ref{eq3.1.6}) that $M$, $N$ and
$[M,N]$ are linearly independent if and only if $[M,N]$ is not
parabolic. ~~~~~~~~~~~~~~~$\square$

\noindent\begin{lemma}\label{lemma3.1.3}%
Assume that $M$ and $N$ are linearly independent elements of
$su(1,1)$ and the set
$\{u\in\mathbb{R}|\left<M\right.+uN,M^{\dag}+\left.uN^{\dag}\right><0\}$
is empty, then $M+uN$ is hyperbolic for each $u\in\mathbb{R}$ if
the commutator $[M,N]$ is not parabolic.
\end{lemma}

{\it{Proof:~}} Because the set
$\left\{u\in\mathbb{R}|\left<M+uN,M^{\dag}+uN^{\dag}\right><0\right\}$
is empty, we have%
\noindent\begin{equation}\label{eq3.1.9}
\left<M+uN,M^{\dag}+uN^{\dag}\right>\geq0,~\forall~u\in\mathbb{R},
\end{equation}

\noindent or equivalently%
\noindent\begin{equation}\label{eq3.1.10}
\left<N,N^{\dag}\right>u^2+2\left<M,N^{\dag}\right>u+\left<M,M^{\dag}\right>\geq0,~\forall~u\in\mathbb{R}.
\end{equation}

\noindent The case that $\left<N,N^{\dag}\right><0$ can be
directly excluded from (\ref{eq3.1.10}). For the case when
$\left<N,N^{\dag}\right>=0$, from (\ref{eq3.1.10}) we have
$\left<M,N^{\dag}\right>=0$. Then, combined with (\ref{eq3.1.6}),
$[M,N]$ must be parabolic, which contradicts with the assumption.
For the case when $\left<N,N^{\dag}\right>>0$, (\ref{eq3.1.10})
holds if and only if
$\left<M,N^{\dag}\right>^2-\left<N,N^{\dag}\right>\left<M,M^{\dag}\right>\leq0$.
If $[M,N]$ is not parabolic, the previous inequality can be
rewritten as
$\left<M,N^{\dag}\right>^2-\left<N,N^{\dag}\right>\left<M,M^{\dag}\right><0$,
which implies that $M+uN$ is hyperbolic for each $u\in\mathbb{R}$.
~~~~~~~~~~~~~~~$\square$

\subsection{Controllability for Single-Input
Case}\label{sec3.2}

Assume that there is only one control in (\ref{eq1.1}), i.e.,%
\noindent\begin{equation}\label{eq3.2.1}
\dot{X}(t)=[A+u(t)B]X(t),~X(0)=I_2.
\end{equation}

\noindent If $A$ and $B$ are linearly independent, i.e., they
commute with each other, the solution of system (\ref{eq3.2.1})
can be expressed as%
\noindent\begin{equation}\label{eq3.2.2}
X(t)=\exp\left[At+B\int_0^tu(\tau)d\tau\right].
\end{equation}

\noindent Accordingly, the reachable set
$R(\cup{\infty}){\subseteq}e^\mathcal{L}=\{X|X=\exp(Bs),s\in\mathbb{R}\}$
is a proper subgroup of $SU(1,1)$, which can never fill up
$SU(1,1)$. Thus, system (\ref{eq3.2.1}) is always uncontrollable
in this case. In the following, we only consider the nontrivial
case when $A$ and $B$ are linearly independent.

For systems evolving on the compact Lie group $SU(2)$, it has been
shown in \cite{Alessandro1, Ramakrishna2} that linear independence
of $A$ and $B$ is a sufficient condition for the involved system
to be controllable. But for the noncompact case of $SU(1,1)$, the
situation is much more complicated. In fact, we have:

\noindent\begin{theorem}\label{theorem3.2.1}%
System (\ref{eq3.2.1}) is uncontrollable if $[A,B]$ is parabolic.
\end{theorem}

{\it{Proof:~}} According to Lemma~\ref{lemma3.1.2}, $A$, $B$ and
$[A,B]$ are linearly dependent when $[A,B]$ is parabolic, which
implies that the Lie algebra
$\mathcal{L}=\{A,B\}_{LA}=\text{span}\{A,B\}$ is two dimensional
and never fills up $su(1,1)$. Thus, the system (\ref{eq3.2.1}) is
uncontrollable on $SU(1,1)$ when $[A,B]$ is parabolic.
~~~~~~~~~~~~~~~$\square$

In addition, even when $[A,B]$ is not parabolic which means that
$A$ and $B$ can generate the whole Lie algebra $su(1,1)$, the
system (\ref{eq3.2.1}) still may be uncontrollable
.%
\noindent\begin{theorem}\label{theorem3.2.2}%
Assume that $[A,B]$ is not parabolic, the system (\ref{eq3.2.1})
is uncontrollable if $A+uB$ is hyperbolic for all
$u\in\mathbb{R}$.
\end{theorem}

{\it{Proof:~}} Since $A+uB$ is hyperbolic for each
$u\in\mathbb{R}$, we have%
\noindent\begin{equation}\label{eq3.2.3}%
\left<B,B^\dag\right>u^2+2\left<A,B^\dag\right>u+\left<A,A^\dag\right>>0,~{\forall}u\in\mathbb{R}.
\end{equation}

\noindent Since $[A,B]$ is not parabolic, from (\ref{eq3.2.3}), we
can immediately obtain that $\left<A,A^\dag\right>>0$ and
$\left<B,B^\dag\right>>0$. Since $B$ is hyperbolic, $B$ can be
converted into $\sqrt{\left<B,B^\dag\right>}\bar{K}_y$ through a
transformation $P$ selected from $SU(1,1)$ (See the Appendix for
rigorous proof). This induces a coordinate transformation in
$SU(1,1)$, given by $X\rightarrow{P^{-1}XP}$, under which the
system (\ref{eq3.2.1}) can be changed into%
\noindent\begin{equation}\label{eq3.2.4}
\dot{\tilde{X}}(t)=[\tilde{A}+\tilde{u}(t)\bar{K}_y]\tilde{X}(t),~\tilde{X}(0)=I_2,
\end{equation}

\noindent where $\tilde{X}=P^{-1}XP$, $\tilde{A}=P^{-1}AP$ and
$\tilde{u}=\sqrt{\left<B,B^\dag\right>}u$. Without loss of
generality, it can be assumed that
$\left<\tilde{A},\bar{K}_y^\dag\right>=0$. In fact, if
$\left<\tilde{A},\bar{K}_y^\dag\right>\neq0$, we can write $u(t)$
in (\ref{eq3.2.2}) as
$u(t)=v(t)-\left<\tilde{A},\bar{K}_y^\dag\right>$ and regard
$\tilde{A}-\left<\tilde{A},\bar{K}_y^\dag\right>\bar{K}_y$ as the
new drift term and $v(t)$ as the new control function. Thus, we
can express $\tilde{A}$ as $a_x\bar{K}_x+a_z\bar{K}_z$, where
$|a_x|>|a_z|$ because $A$ is hyperbolic. Rescaling the time
variable $t$ by a factor $|a_x|$ gives a system of the form as%
\noindent\begin{equation}\label{eq3.2.5}
\dot{{X}}(t)=[\varepsilon\bar{K}_x+a\bar{K}_z+u(t)\bar{K}_y]X(t),~X(0)=I_2,
\end{equation}

\noindent where $\varepsilon=\text{sgn}(a_x)=\pm1$ and $|a|<1$ .
Clearly, system (\ref{eq3.2.5}) shares the same controllability
properties with system (\ref{eq3.2.1}).

Now, we prove that system (\ref{eq3.2.5}) is uncontrollable. Write
the solution of the evolution equation (\ref{eq3.2.5}) as%
\noindent\begin{equation}\label{eq3.2.6}
X:=\left(%
\begin{array}{cc}
  x_1+ix_2 & x_3-ix_4 \\
  x_3+ix_4 & x_1-ix_2
\end{array}%
\right),
\end{equation}

\noindent then we have%
\noindent\begin{equation}\label{eq3.2.7}
  \dot{x}_1=\frac{1}{2}(ax_2+{\varepsilon}x_4-ux_3),
\end{equation}
\noindent\begin{equation}\label{eq3.2.8}
  \dot{x}_2=\frac{1}{2}(-ax_1-{\varepsilon}x_3-ux_4),
\end{equation}
\noindent\begin{equation}\label{eq3.2.9}
  \dot{x}_3=\frac{1}{2}(-ax_4-{\varepsilon}x_2-ux_1),
\end{equation}
\noindent\begin{equation}\label{eq3.2.10}
  \dot{x}_4=\frac{1}{2}(ax_3+{\varepsilon}x_1-ux_2).
\end{equation}

\noindent Subtracting Eqs.(\ref{eq3.2.7}) and (\ref{eq3.2.10})
then followed by a succeeding multiplication by $2(x_1-x_4)$ gives%
\noindent\begin{equation}\label{eq3.2.11}
  \frac{d}{dt}(x_1-x_4)^2=a(x_1-x_4)(x_2-x_3)-{\varepsilon}(x_1-x_4)^2+u(x_1-x_4)(x_2-x_3).
\end{equation}

\noindent Similarly, we have%
\noindent\begin{equation}\label{eq3.2.12}
  \frac{d}{dt}(x_2-x_3)^2={-a}(x_1-x_4)(x_2-x_3)+{\varepsilon}(x_2-x_3)^2+u(x_1-x_4)(x_2-x_3).
\end{equation}

\noindent Then, subtracting Eqs.(\ref{eq3.2.11}) and
(\ref{eq3.2.12}) derives%
\noindent\begin{equation}\label{eq3.2.13}
\begin{split}
  \frac{d}{dt}[(x_1-x_4)^2-(x_2-x_3)^2]&=2a(x_1-x_4)(x_2-x_3)-{\varepsilon}[(x_1-x_4)^2+(x_2-x_3)^2] \\
  &=-{\varepsilon}(1-|a|)[(x_1-x_4)^2+(x_2-x_3)^2]-{\varepsilon}|a|[(x_1-x_4)-\text{sgn}(a){\varepsilon}(x_2-x_3)]^2\\
  &\left\{%
    \begin{array}{ll}
        \leq0, & \hbox{when $\varepsilon=1$;} \\
       \geq0, & \hbox{when $\varepsilon=-1$.} \\
    \end{array}%
\right.
\end{split}
\end{equation}

\noindent Thus the function $[x_1(t)-x_4(t)]^2-[x_2(t)-x_3(t)]^2$
is nonincreasing (nondecreasing) for every trajectory of system
(\ref{eq3.2.5}) when $\varepsilon=1$ ($\varepsilon=-1$). Since the
initial value of this function is $1$, it can be concluded that
the reachable states of system (\ref{eq3.2.5}) should satisfy the
restriction $(x_1-x_4)^2-(x_2-x_3)^2\leq1 (\geq1)$ when
$\varepsilon=1$ ($\varepsilon=-1$). This result means that the
reachable set of system (\ref{eq3.2.5}) never equals $SU(1,1)$,
i.e., the involved system is uncontrollable. This completes the
proof. ~~~~~~~~~~~~~~~$\square$

Combining Lemma~\ref{lemma3.1.3} and Theorems~\ref{theorem3.2.1}
and \ref{theorem3.2.2}, we can immediately obtain the following
result.

\noindent\begin{theorem}\label{theorem3.2.3}%
System (\ref{eq3.2.1}) is uncontrollable if the set
\end{theorem}

\noindent\begin{equation}\label{eq3.2.14}
\Omega=\left\{u\in\mathbb{R}|\left<A+uB,A^{\dag}+uB^{\dag}\right><0\right\}
\end{equation}

\noindent is empty.

This theorem suggests that only when the operator $A+uB$ can be
adjusted to be elliptic by some constant $u\in\mathbb{R}$ can we
realize arbitrary propagators of the system as an element in the
noncompact Lie group $SU(1,1)$.

When the admissible control set $\mathcal{U}$ is assumed to be the
class of all locally bounded and measurable functions, a
sufficient condition is given in \cite{Jurdjevic1} for the
controllability of the system on more general Lie groups. This
condition states that the involved system is controllable if there
exists a constant control $u$ such that the resulting state
trajectory is periodic in the course of time. Since $\exp(tM)$ is
periodic if and only if it is elliptic, we can extend this result
to the case of $SU(1,1)$ as follows.

\noindent\begin{theorem}\label{theorem3.2.4}%
System (\ref{eq3.2.1}) is controllable if and only if the set
$\Omega$ in (\ref{eq3.2.14}) is nonempty.
\end{theorem}

This theorem means that the controllability system (\ref{eq3.2.1})
is completely characterized by the set $\Omega$, and thus provides
a sufficient and necessary condition that examines the
controllability of single-input control system on $SU(1,1)$. Since
the value of the set $\Omega$ is completely determined by $A$ and
$B$, we can further describe the system controllability with
respect to $A$ and $B$ as specified in the following table.

\noindent\begin{center}
\begin{tabular}{c|c|c}
\multicolumn{3}{c}{\text{Table~I. The controllability
characterization of
system~(\ref{eq3.2.1}).}}\\
  \hline\hline
  \rule[0pt]{0cm}{5pt}%
   The range of $A$ and $B$ & The set $\Omega$ & System controllability \\\hline
$\begin{array}{c}%
  \left<B,B^\dag\right><0,\\
  \left<A,B^\dag\right>^2-\left<A,A^\dag\right>\left<B,B^\dag\right>\neq0
\end{array}$%
    &~~Nonempty~~& Controllable \\\hline
  $\left<B,B^\dag\right>=0$,~~~$\left<A,B^\dag\right>\neq0$ &~~Nonempty~~& Controllable \\\hline
$\begin{array}{c}%
  \left<B,B^\dag\right>>0,\\
  \left<A,B^\dag\right>^2-\left<A,A^\dag\right>\left<B,B^\dag\right>>0
\end{array}$%
    &~~Nonempty~~& Controllable \\\hline
    Otherwise &   Empty    & Uncontrollable \\\hline\hline
\end{tabular}
\end{center}

\noindent{\bf{Remark:}} If the admissible control $u(t)$ is
restricted by an up-bound, i.e., $|u(t)|{\leq}C$ for any $t\geq0$,
where $C$ is a priori prescribed positive constant, a similar
conclusion can be drawn for the system (\ref{eq3.2.1}). The
relevant necessary and sufficient condition can be constructed
by the following set%
\noindent\begin{equation}\label{eq3.2.16}
\tilde{\Omega}=\left\{-C{\leq}u{\leq}C\left|\left<A+uB,A^{\dag}+uB^{\dag}\right><0\right\}\right..
\end{equation}

\noindent It was shown in \cite{wjw1} that any element
$X_f{\in}SU(1,1)$ can be decomposed as%
\noindent\begin{equation}\label{eq3.2.16}
X_f=\prod\limits_{k=1}^Q\exp[T_k(A+u_kB)]
\end{equation}

\noindent  when $\tilde{\Omega}$ is nonempty, where $T_k\geq0$,
$u_k{\leq}C$ and $Q$ is a positive integer number. This result
indicates that the nonemptiness of the set $\tilde{\Omega}$ is the
corresponding sufficient condition for the controllability of the
system. This condition also can be proved to be necessary in a
similar way as that of Theorem~\ref{theorem3.2.2} (see Example~2
for illustration).

Now, we turn to the strong controllability. In the following, we
will show that system (\ref{eq3.2.1}) is never strong
controllable. Without loss of generality, we assume that the
admissible controls are piecewise constant functions of $t$ with a
finite number of switches, i.e., any time interval $[0, t_f]$ can
be partitioned into $N$ subintervals $[t_{k-1}, t_k]$ such that
$t_0=0$, $t_N=t_f$ and any control $u(t)$ takes a constant value
$u_k$ on $(t_{k-1}, t_k)$. Accordingly, the time evolution of
system (\ref{eq3.2.1}) can be
expressed as%
\noindent\begin{equation}\label{eq3.2.17}
X(u(\cdot),t_f)=\prod_{k=1}^{N}\exp[T_k(A+u_kB)],
\end{equation}

\noindent where $T_k=t_k-t_{k-1}$. Since%
\noindent\begin{equation}\label{eq3.2.18}%
\lim\limits_{T_k\rightarrow0}e^{T_k(A+u_kB)}\left\{%
\begin{array}{ll}
    =I_2, & \hbox{when $\lim\limits_{T_k\rightarrow0}u_kT_k=0$;} \\
    \in{\left\{e^{sB}|s\neq0\right\}}, & \hbox{otherwise,} \\
\end{array}%
\right.
\end{equation}

\noindent we have, for any given $u(t)$,
\noindent\begin{equation}\label{eq3.2.18}%
\lim\limits_{t_f\rightarrow0}X(u(\cdot),t_f)\in{\left\{e^{sB}|s\in\mathbb{R}\right\}}.
\end{equation}

\noindent Thus,
$R(\cap{\infty})\subseteq{\left\{e^{sB}|s\in\mathbb{R}\right\}}$,
i.e., system (\ref{eq3.2.1}) is not strong controllable.

Since for any given time $t_f$ and $s\in\mathbb{R}$, we can choose
a constant control $\bar{u}=\frac{s}{t_f}$, and then have
$\lim\limits_{t_f\rightarrow0}e^{t_f(A+\bar{u}B)}=e^{sB}$. Thus,
we have
${\left\{e^{sB}|s\in\mathbb{R}\right\}}\subseteq\lim\limits_{t\rightarrow0}\overline{R(t)}$,
and can further draw the conclusion that
${\left\{e^{sB}|s\in\mathbb{R}\right\}}{\subseteq}\overline{R(\cap{\infty})}$
when system (\ref{eq3.2.1}) is small time local controllable.

We have the following result for the small time local
controllability of system (\ref{eq3.2.1}).

\noindent\begin{theorem}\label{theorem3.2.5}%
System (\ref{eq3.2.1}) is small time local controllable if
$\left<B,B^\dag\right><0$.
\end{theorem}

{\it{Proof:~}} Since $\left<B,B^\dag\right><0$, there exists a
positive quantity $u_c$ such that
$\left<A+uB,A^\dag+uB^\dag\right><0$ for every $u>u_c$. When
$u>u_c$, the eigenvalues of $(A+uB)\varepsilon$ are
$\lambda_{1,2}={\pm}{i\varepsilon}\sqrt{-\frac{1}{4}\left(\left<A,A^\dag\right>+2\left<A,B^\dag\right>u+\left<B,B^\dag\right>u^2\right)}$
for each $\varepsilon>0$. Thus, the value of $u$ can be chosen
such that $\lambda_{1,2}={\pm}i2n\pi$, so we have
$e^{\varepsilon(A+uB)}=I_2$. Since $u$ is nonzero, it can be
proved that $I_2$ is an interior point of $R(\varepsilon)$ with
the similar method used in~\cite{Alessandro1}. Thus system
(\ref{eq3.2.1}) is small time local controllable if
$\left<B,B^\dag\right>$ is negative. ~~~~~~~~~~~~~~~$\square$

\subsection{Controllability for Multi-Input
Case}\label{sec3.3}

In this section, we consider the controllability of system
(\ref{eq1.1}) with multiple inputs. Since the matrices
$B_k,k=1,\cdots,r$, have been assumed to be linearly independent,
it is sufficient to consider the following two cases: (I) $r=3$,
it is obvious that $B_1$, $B_2$ and $B_3$ generate the whole Lie
algebra of $su(1,1)$, and we have
$\mathcal{L}=\mathcal{L}_0=\mathcal{B}=su(1,1)$, which means that
system (\ref{eq1.1}) is strong controllable; (II) $r=2$, i.e.,
\noindent\begin{equation}\label{eq3.3.1}
\dot{X}(t)=[A+u_1(t)B_1+u_2(t)B_2]X(t),~X(0)=I_2,
\end{equation}

\noindent for which we have
\noindent\begin{theorem}\label{theorem8.1}%
\noindent\begin{itemize}
    \item[i)] If $A$ can be written as linear combination of $B_1$ and
$B_2$, then system (\ref{eq3.3.1}) is uncontrollable if
$[B_1,B_2]$ is parabolic. Otherwise, it is strong controllable.
    \item[ii)] If $A$, $B_1$ and $B_2$ are linearly independent, then
system (\ref{eq3.3.1}) is controllable. Moreover, it is strong
controllable if $[B_1,B_2]$ is not parabolic.
\end{itemize}
\end{theorem}%

{\it{Proof:~}} i) Since $A$ can be written as linear combination
of $B_1$ and $B_2$, according to Lemma~\ref{lemma3.1.2}, $A$,
$B_1$ and $B_2$ do not generate the whole Lie algebra of $su(1,1)$
when $[B_1,B_2]$ is parabolic, i.e.,
$\mathcal{L}=\left\{A,B_1,B_2\right\}_{LA}=\left\{B_1,B_2\right\}_{LA}=\text{span}\left\{B_1,B_2\right\}{\neq}su(1,1)$.
This means that system (\ref{eq3.3.1}) is not controllable. When
$[B_1,B_2]$ is not parabolic, since accordingly $B_1$, $B_2$ and
$[B_1,B_2]$ form a basis in $su(1,1)$, we have
$\mathcal{B}=\left\{B_1,B_2\right\}_{LA}=su(1,1)$. This implies
that system (\ref{eq3.3.1}) is strong controllable. ii) Since $A$,
$B_1$ and $B_2$ are linearly independent, there must exist two
constants $\bar{u}_1$ and $\bar{u}_2$ such that
$A+\bar{u}_1B_1+\bar{u}_2B_2$ is elliptic. Thus, from the results
obtained in the previous subsection, we can conclude that system
(\ref{eq3.3.1}) is controllable. A similar argument as in i) can
be given to the case that $[B_1,B_2]$ is not parabolic to show
that system (\ref{eq3.3.1}) is strong controllable.
~~~~~~~~~~~~~~~$\square$

\section{Relation Between Systems on $SU(1,1)$, $SO(2,1)$ and $SL(2,\mathbb{R})$}\label{sec4}

In this section, we show that the results obtained in
Section~\ref{sec3} are also valid for the systems on the Lie
groups $SO(2,1)$ and $SL(2,\mathbb{R})$, because both the map
$\rho_1:~su(1,1)~\rightarrow~so(2,1)$ defined by%
\noindent\begin{equation}\label{eq4.1}%
\rho_1:=\bar{K}_\alpha~\rightarrow~O_\alpha,~\alpha=x,y,z,
\end{equation}

\noindent with%
\noindent\begin{equation}\label{eq4.2}%
O_x=\left(\begin{array}{ccc}
  0 & 0 & 0 \\
  0 & 0 & 1 \\
  0 & 1 & 0 \\
\end{array}\right),~%
O_y=\left(\begin{array}{ccc}
  0 & 0 & 1 \\
  0 & 0 & 0 \\
  1 & 0 & 0 \\
\end{array}\right),~%
O_z=\left(\begin{array}{ccc}
  0 & 1 & 0 \\
  -1& 0 & 0 \\
  0 & 0 & 0 \\
\end{array}\right),
\end{equation}

\noindent and the map
$\rho_2:~su(1,1)~\rightarrow~sl(2,\mathbb{R})$ defined by%
\noindent\begin{equation}\label{eq4.3}%
\rho_2:=\bar{K}_\alpha~\rightarrow~L_\alpha,~\alpha=x,y,z,
\end{equation}

\noindent with%
\noindent\begin{equation}\label{eq4.4}%
L_x=\frac{1}{2}\left(\begin{array}{cc}
  -1 & 0 \\
  0 & 1 \\
\end{array}\right),~%
L_y=\frac{1}{2}\left(\begin{array}{cc}
  0 & -1 \\
  -1 & 0\\
\end{array}\right),~%
L_z=\frac{1}{2}\left(\begin{array}{cc}
  0 & -1 \\
  1 & 0 \\
\end{array}\right),%
\end{equation}

\noindent are Lie algebra isomorphism. According to Lie's third
theorem, $\rho_1$ and $\rho_2$ induce a two-to-one homomorphism
$\tilde{\rho}_1$ from $SU(1,1)$ to $SO(2,1)$ and a isomorphism
$\tilde{\rho}_2$ from $SU(1,1)$ to $SL(2,\mathbb{R})$
respectively~\cite{Vilenkin1}. Accordingly, we can associate the
system given in (\ref{eq1.1}) to
the system varying on $SO(2,1)$%
\noindent\begin{equation}\label{eq4.5}
\dot{Y}(t)=[\rho_1(A)+\sum_{l=1}^{r}u_l(t)\rho_1(B_l)]Y(t),~Y(0)=I_3,
\end{equation}

\noindent and the system varying on $SL(2,\mathbb{R})$%
\noindent\begin{equation}\label{eq4.6}
\dot{Z}(t)=[\rho_2(A)+\sum_{l=1}^{r}u_l(t)\rho_2(B_l)]Z(t),~Z(0)=I_2,
\end{equation}

\noindent respectively. The state of system (\ref{eq4.5}) consists
of all the transformations that leave the three-dimensional
hyperboloids $x^2+y^2-z^2=\pm1$ invariant, while the state of
system (\ref{eq4.6}) consists of all the $2\times2$ real matrices
with determinant 1. Clearly, when we impose the same controls
$u_l(t)$ on the systems (\ref{eq1.1}), (\ref{eq4.5}) and
(\ref{eq4.6}), their trajectories can be mapped by
$\tilde{\rho}_1$ and $\tilde{\rho}_2$ respectively, i.e.,
$Y(t;u_l(\cdot))=\rho_1(X(t;u_l(\cdot)))$ and
$Z(t;u_l(\cdot))=\rho_2(X(t;u_l(\cdot)))$. Therefore, the
controllability properties of the associated systems (\ref{eq4.5})
and (\ref{eq4.6}) can be obtained from system (\ref{eq1.1})
directly.

This also provides a way of picturing the control over Lie group
$SU(1,1)$ by project it onto $SO(2,1)$ as shown in Fig.\ref{fig1}.
The problem of steering system (\ref{eq1.1}) to an arbitrary state
$X_f$ from the initial state $I_2$ can be viewed as the problem of
finding a path between two arbitrary points $P_1$ and $P_2$ on the
hyperboloid of one sheet. As shown in Fig.\ref{fig1}, the
$SO(2,1)$ evolution operators $e^{tO_{\alpha}}$ ($\alpha=x,y,z$)
are identified with the rotations about $\alpha$-axis. Thus,
piecewise constant controls induce a series of rotations about the
axis through the origin $O$. For example, when system
(\ref{eq3.2.1}) is under the action of constant control $u$, the
induced rotation is $e^{t[\rho_1(A)+u\rho_1(B)]}$. Because the
evolution time is assigned to be nonnegative, the rotation induced
can be performed only in one direction. Theorem~\ref{theorem3.2.4}
suggests that, if and only if the system can rotate about at least
one axis that is located inside the cone $x^2+y^2-z^2\leq0$, can
we move any given point on the hyperboloid to another one via a
series of rotations. Under the rotation about the axis that is
located inside the cone, every point on the hyperboloid follows a
closed
elliptic trajectory.%
\noindent\begin{center} \noindent\begin{figure}[h] \centering
\includegraphics[totalheight=4.0in]{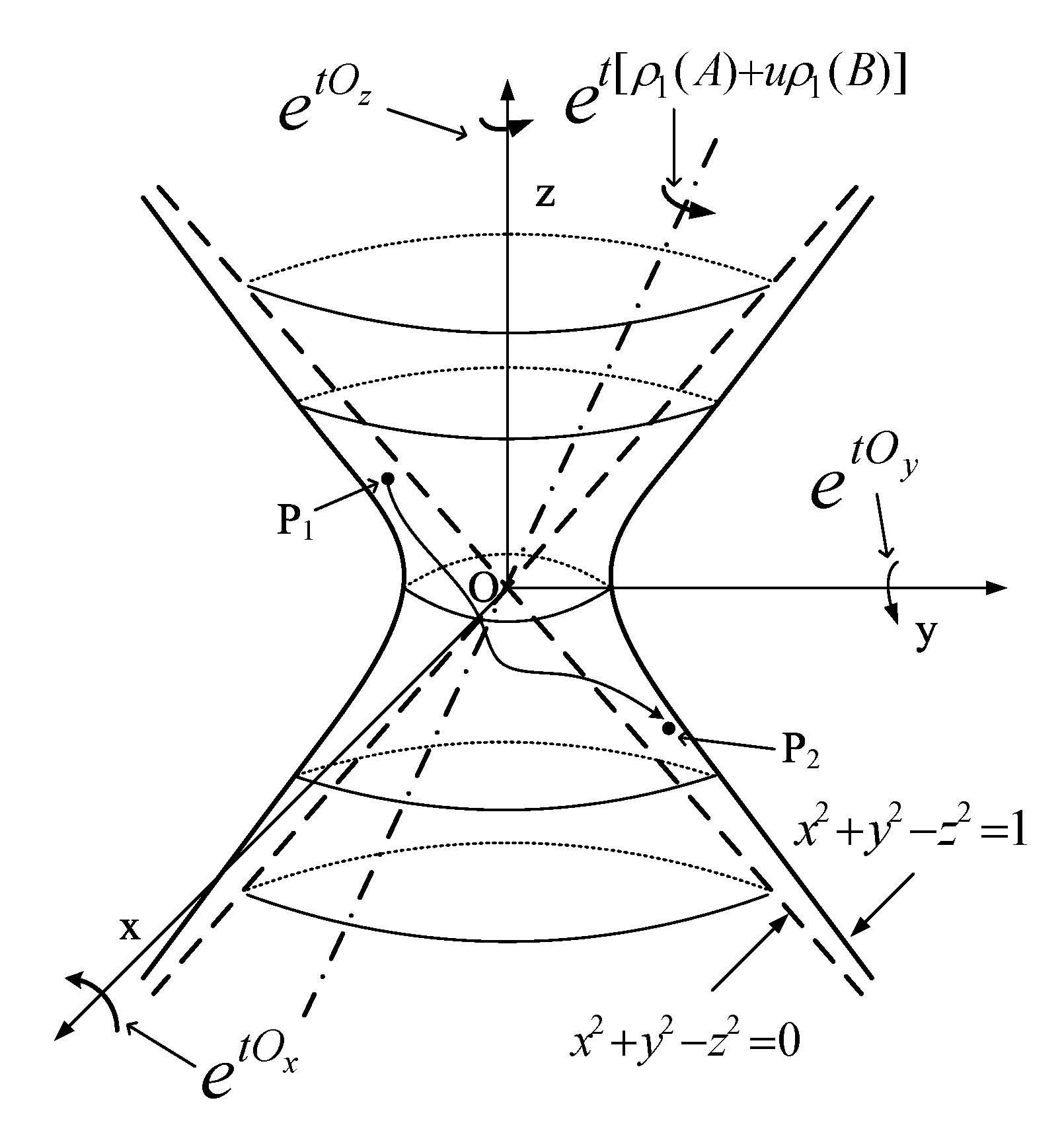}%{fig.eps}%
\caption{The topology of $SO(2,1)$.}\label{fig1}
\end{figure}
\end{center}\noindent

\section{Examples}\label{sec5}

\noindent{\bf{Example~1:}} Consider the quantum system with its
Hamiltonian expressed as~\cite{Gerry2}%
\noindent\begin{equation}\label{eq5.1}%
H(t)=\omega_0K_z+u(t)K_x,
\end{equation}

\noindent where $K_x$ and $K_z$ are operators as defined in
(\ref{eq2.1.2}). The quantum system%
\noindent\begin{equation}\label{eq5.2}%
i\hbar\dot{\psi}(t)=H(t)\psi(t)
\end{equation}

\noindent is then a quantum control system that preserves
$SU(1,1)$ coherent states~\cite{Gerry1}. Consider the positive
discrete series unitary irreducible representations of $su(1,1)$
denoted by $\mathscr{D}^{+}(k)$, where $k$ is the so-called
Bargmann index. The basis states $|m,k\rangle$ diagonalize the
generator $K_z$ and the Casimir operator $C=K_z^2-K_x^2-K_y^2$ as
follows: $K_z|m,k\rangle=(m+k)|m,k\rangle~(m=0,1,2,\cdots)$, and
$C|m,k\rangle=k(k-1)|m,k\rangle$ with $k>0$. Then the operators
$K_{\pm}=K_x{\pm}iK_y$ will act as raising and lowering operators,%
\noindent\begin{equation}\label{eq5.3}
\begin{array}{l}
   K_+|m,k\rangle=[(m+1)(m+2k)]^{1/2}|m+1,k\rangle,\\
   K_-|m,k\rangle=[m(m+2k-1)]^{1/2}|m-1,k\rangle.\\
\end{array}
\end{equation}

\noindent With the representation introduced above, the operators
$K_\pm$ and $K_z$ are then identified as
\noindent\begin{equation}\nonumber%
K_+=\left(%
\begin{array}{ccccc}
  0 &  \sqrt{2k} &  &  &  \\
    & 0 & 2\sqrt{2k+1} &  &  \\
    &   & 0 & 3\sqrt{2k+2} &  \\
    &   &  & 0 & \ddots \\
    &   &  &  & \ddots \\
\end{array}%
\right),
\end{equation}
\noindent\begin{equation}\nonumber%
K_-=\left(%
\begin{array}{ccccc}
  0 &   &  &  &  \\
  \sqrt{2k} & 0 &  &  &  \\
    &  2\sqrt{2k+1} & 0 &  &  \\
    &   & 3\sqrt{2k+2} & 0 &  \\
    &   &  & \ddots & \ddots \\
\end{array}%
\right),
\end{equation}
\noindent\begin{equation}\nonumber%
K_z=\left(%
\begin{array}{ccccc}
  k+1 &  &  &  &  \\
   & k+2 &  &  &  \\
   &  & k+3 &  &  \\
   &  &  & k+4 &  \\
   &  &  &  & \ddots\\
\end{array}%
\right).
\end{equation}

\noindent Following Perelomov~\cite{Perelomov1}, the $SU(1,1)$
coherent states are expressed as a linear combination of the basis
vectors $|m,k\rangle~(m=0,1,2,\cdots)$, and can be obtained from
the state $|0,k\rangle$ by the action of
$\exp(\alpha{K_+}-\alpha^*{K_-})=\exp\{-2[\text{Im}(\alpha)(-iK_x)+\text{Re}(\alpha)(-iK_y)]\}$,
where $\alpha$ is a complex number. Since, according to
Theorem~\ref{theorem3.2.4}, the equivalent system of the system
(\ref{eq5.2})%
\noindent\begin{equation}\label{eq5.4}%
\dot{X}(t)=[\omega_0\bar{K}_z+u(t)\bar{K}_x]X(t)
\end{equation}

\noindent is controllable on $SU(1,1)$, it can be concluded that
the transition between two arbitrary $SU(1,1)$ coherent states can
be realized by controlling the quantum system (\ref{eq5.2}).

\noindent{\bf{Example~2:}} Consider the following control system
evolving on $SO(2,1)$~\cite{Jurdjevic2}%
\noindent\begin{equation}\label{eq5.5}%
\dot{Y}(t)=[O_x+u(t)O_z]Y(t),~Y(0)=I_3,%
\end{equation}

\noindent and assume that the control $u(t)$ is restricted by
$|u(t)|{\leq}C$, then the system is controllable if and only if
$C>1$.

The associated system, evolving on $SU(1,1)$, is as follows
\noindent\begin{equation}\label{eq5.6}%
\dot{X}(t)=[\bar{K}_x+u(t)\bar{K}_z]X(t),~X(0)=I_2.%
\end{equation}

\noindent It can be verified that the set
$\bar{\Omega}=\left\{-C{\leq}u{\leq}C|\left<\bar{K}_x+u\bar{K}_z,\bar{K}_x^\dag+u\bar{K}_z^\dag\right>\right\}$
is nonempty if and only if $C>1$. Thus, according to the results
obtained in Section~\ref{sec3}, system (\ref{eq5.6}) is
controllable when $C>1$.

Now we show that system (\ref{eq5.6}) is uncontrollable when
$C{\leq}1$. Write the solution of the evolution equation
(\ref{eq5.6}) as%
\noindent\begin{equation}\label{eq5.7}
X:=\left(%
\begin{array}{cc}
  x_1+ix_2 & x_3-ix_4 \\
  x_3+ix_4 & x_1-ix_2
\end{array}%
\right),
\end{equation}

\noindent then, with a few calculations, we have%
\noindent\begin{equation}\label{eq5.8}
\begin{split}
  \frac{d}{dt}[(x_1-x_4)^2-(x_2-x_3)^2]&=2u(x_1-x_4)(x_2-x_3)-[(x_1-x_4)^2+(x_2-x_3)^2] \\
  &=-(1-u^2)(x_1-x_4)^2-[u(x_1-x_4)-(x_2-x_3)]^2\\
  &\leq0.
\end{split}
\end{equation}

\noindent This means that the function
$[x_1(t)-x_4(t)]^2-[x_2(t)-x_3(t)]^2$ is nonincreasing for every
trajectory of system (\ref{eq5.6}) if $|u|{\leq}1$. Thus, the
reachable set of system (\ref{eq5.6}) never equals $SU(1,1)$, and
the system is accordingly uncontrollable. As a result, system
(\ref{eq5.5}) is controllable if and only if $C>1$.

\section{Conclusion}\label{sec6}

In this paper, we have studied the controllability properties of
the quantum system evolving on the noncompact Lie group $SU(1,1)$.
The criteria established in this article can be used to examine,
for example, the ability to control the transitions between
different $SU(1,1)$ coherent states. The results obtained in this
paper also can be extended to the systems evolving on $SO(2,1)$
and $SL(2,\mathbb{R})$, because they are both homomorphic to
$SU(1,1)$.

\section*{Acknowledgment}

The authors would like to thank Dr. Re-Bing Wu for his helpful
suggestions.

\appendix
\section*{Appendix}
In this appendix, we show that any hyperbolic $B$ can be converted
into $\sqrt{\left<B,B^\dag\right>}\bar{K}_y$ through a matrix
$P{\in}SU(1,1)$, i.e.,
$PBP^{-1}=\sqrt{\left<B,B^\dag\right>}\bar{K}_y$. Since $B$ is
hyperbolic, we can expand it in the basis given in (\ref{eq2.1.4})
as $B=x\bar{K}_x+y\bar{K}_y+z\bar{K}_z$, where
$\left<B,B^\dag\right>=x^2+y^2-z^2>0$.

First, one can find a matrix
$P_1=e^{\alpha\bar{K}_z}{\in}SU(1,1)$, which satisfy%
\noindent\begin{equation}\label{eqa.1}%
P_1BP_1^{-1}=\sqrt{x^2+y^2}\bar{K}_y+z\bar{K}_z.
\end{equation}

\noindent Let $\alpha$ be the angle satisfying%
\noindent\begin{equation}\label{eqa.2}%
\sin\alpha=\frac{x}{\sqrt{x^2+y^2}},~~\cos\alpha=\frac{y}{\sqrt{x^2+y^2}}.
\end{equation}

\noindent According to the Baker-Hausdorff-Campbell formula%
\noindent\begin{equation}\label{eqa.3}%
e^{M}Ne^{-M}=N+[M,N]+\frac{1}{2!}[M,[M,N]]+\frac{1}{3!}[M,[M,[M,N]]]+\cdots,
\end{equation}

\noindent one can immediately obtain that
\noindent\begin{equation}\label{eqa.4}%
\begin{array}{l}
   e^{\alpha\bar{K}_z}Be^{-\alpha\bar{K}_z}=xe^{\alpha\bar{K}_z}\bar{K}_xe^{-\alpha\bar{K}_z}+ye^{\alpha\bar{K}_z}\bar{K}_ye^{-\alpha\bar{K}_z}+z\bar{K}_z\\
   ~~~~~~~~~~=(x\cos\alpha-y\sin\alpha)\bar{K}_x+(x\sin\alpha+y\cos\alpha)\bar{K}_y+z\bar{K}_z\\
   ~~~~~~~~~~=\sqrt{x^2+y^2}\bar{K}_y+z\bar{K}_z.
\end{array}
\end{equation}

Next, we show that there is a matrix $P_2=e^{\beta\bar{K}_x}$, in
$SU(1,1)$, which can convert $\sqrt{x^2+y^2}\bar{K}_y+z\bar{K}_z$
into $\sqrt{\left<B,B^\dag\right>}\bar{K}_y$. Since
$x^2+y^2-z^2>0$, we can choose $\beta$ such that
\noindent\begin{equation}\label{eqa.5}%
\sinh\beta=\frac{z}{\sqrt{x^2+y^2-z^2}},~~\cosh\beta=\frac{\sqrt{x^2+y^2}}{\sqrt{x^2+y^2-z^2}}.
\end{equation}

\noindent Make use of the formula given in (\ref{eqa.3}) again, we
have
\noindent\begin{equation}\label{eqa.4}%
\begin{array}{l}
   e^{\beta\bar{K}_x}(\sqrt{x^2+y^2}\bar{K}_y+z\bar{K}_z)e^{-\beta\bar{K}_x}\\
   ~~~=\sqrt{x^2+y^2}e^{\beta\bar{K}_x}\bar{K}_ye^{-\beta\bar{K}_x}+ze^{\beta\bar{K}_x}\bar{K}_ze^{-\beta\bar{K}_x}\\
   ~~~=(\sqrt{x^2+y^2}\cosh\beta-z\sinh\beta)\bar{K}_y+(z\cosh\beta-\sqrt{x^2+y^2}\sinh\beta)\bar{K}_z\\
   ~~~=\sqrt{x^2+y^2-z^2}\bar{K}_y\\
   ~~~=\sqrt{\left<B,B^\dag\right>}\bar{K}_y.
\end{array}
\end{equation}

Consequently, the $SU(1,1)$ matrix
$e^{\beta\bar{K}_x}e^{\alpha\bar{K}_z}$ will convert $B$ into
$\sqrt{\left<B,B^\dag\right>}\bar{K}_y$ when it is hyperbolic.


\begin{thebibliography}{23}
\expandafter\ifx\csname
natexlab\endcsname\relax\def\natexlab#1{#1}\fi
\expandafter\ifx\csname bibnamefont\endcsname\relax
  \def\bibnamefont#1{#1}\fi
\expandafter\ifx\csname bibfnamefont\endcsname\relax
  \def\bibfnamefont#1{#1}\fi
\expandafter\ifx\csname citenamefont\endcsname\relax
  \def\citenamefont#1{#1}\fi
\expandafter\ifx\csname url\endcsname\relax
  \def\url#1{\texttt{#1}}\fi
\expandafter\ifx\csname
urlprefix\endcsname\relax\def\urlprefix{URL }\fi
\providecommand{\bibinfo}[2]{#2}
\providecommand{\eprint}[2][]{\url{#2}}

\bibitem[{\citenamefont{Sussmann}(1972)}]{Sussmann1}
\bibinfo{author}{\bibfnamefont{H.~J.} \bibnamefont{Sussmann}},
  \bibinfo{journal}{J. Diff. Eqns.} \textbf{\bibinfo{volume}{12}},
  \bibinfo{pages}{313} (\bibinfo{year}{1972}).

\bibitem[{\citenamefont{R.W.Brockket}(1972)}]{Brockket1}
\bibinfo{author}{\bibnamefont{R.W.Brockket}}, \bibinfo{journal}{SIAM J.
  Control} \textbf{\bibinfo{volume}{10}}, \bibinfo{pages}{265}
  (\bibinfo{year}{1972}).

\bibitem[{\citenamefont{Jurdjevic and Sussman}(1972)}]{Jurdjevic1}
\bibinfo{author}{\bibfnamefont{V.}~\bibnamefont{Jurdjevic}} \bibnamefont{and}
  \bibinfo{author}{\bibfnamefont{H.~J.} \bibnamefont{Sussman}},
  \bibinfo{journal}{J. Diff. Eqns.} \textbf{\bibinfo{volume}{12}},
  \bibinfo{pages}{313} (\bibinfo{year}{1972}).

\bibitem[{\citenamefont{Sussmann}(1983)}]{Sussmann2}
\bibinfo{author}{\bibfnamefont{H.~J.} \bibnamefont{Sussmann}}, in
  \emph{\bibinfo{booktitle}{Differential Geometric Control Theory}}, edited by
  \bibinfo{editor}{\bibfnamefont{R.~W.} \bibnamefont{Brockett}},
  \bibinfo{editor}{\bibfnamefont{R.~S.} \bibnamefont{Millman}},
  \bibnamefont{and} \bibinfo{editor}{\bibfnamefont{H.~J.}
  \bibnamefont{Sussmann}} (\bibinfo{publisher}{Birkhauser},
  \bibinfo{address}{Boston}, \bibinfo{year}{1983}), pp.
  \bibinfo{pages}{1--116}.

\bibitem[{\citenamefont{Huang et~al.}(1983)\citenamefont{Huang, Tarn, and
  Clark}}]{Huang1}
\bibinfo{author}{\bibfnamefont{G.~M.} \bibnamefont{Huang}},
  \bibinfo{author}{\bibfnamefont{T.~J.} \bibnamefont{Tarn}}, \bibnamefont{and}
  \bibinfo{author}{\bibfnamefont{J.~W.} \bibnamefont{Clark}},
  \bibinfo{journal}{J. Math. Phys.} \textbf{\bibinfo{volume}{24}},
  \bibinfo{pages}{2608} (\bibinfo{year}{1983}).

\bibitem[{\citenamefont{Ramakrishna et~al.}(1995)\citenamefont{Ramakrishna,
  Salapaka, Dahleh, Rabitz, and Peirce}}]{Ramakrishna1}
\bibinfo{author}{\bibfnamefont{V.}~\bibnamefont{Ramakrishna}},
  \bibinfo{author}{\bibfnamefont{M.~V.} \bibnamefont{Salapaka}},
  \bibinfo{author}{\bibfnamefont{M.}~\bibnamefont{Dahleh}},
  \bibinfo{author}{\bibfnamefont{H.}~\bibnamefont{Rabitz}}, \bibnamefont{and}
  \bibinfo{author}{\bibfnamefont{A.}~\bibnamefont{Peirce}},
  \bibinfo{journal}{Phys. Rev. A} \textbf{\bibinfo{volume}{51}},
  \bibinfo{pages}{960} (\bibinfo{year}{1995}).

\bibitem[{\citenamefont{Tarn et~al.}(2000)\citenamefont{Tarn, Clark, and
  Lucarelli}}]{Tarn1}
\bibinfo{author}{\bibfnamefont{T.~J.} \bibnamefont{Tarn}},
  \bibinfo{author}{\bibfnamefont{J.~W.} \bibnamefont{Clark}}, \bibnamefont{and}
  \bibinfo{author}{\bibfnamefont{D.~G.} \bibnamefont{Lucarelli}}, in
  \emph{\bibinfo{booktitle}{Proceedings of the 39th IEEE Conference on Decision
  and Control}} (\bibinfo{address}{Sydney}, \bibinfo{year}{2000}), pp.
  \bibinfo{pages}{943--948}.

\bibitem[{\citenamefont{D'Alessandro}(2000)}]{Alessandro1}
\bibinfo{author}{\bibfnamefont{D.}~\bibnamefont{D'Alessandro}},
  \bibinfo{journal}{Systems \& Control Letters} \textbf{\bibinfo{volume}{41}},
  \bibinfo{pages}{213} (\bibinfo{year}{2000}).

\bibitem[{\citenamefont{Schirmer et~al.}(2001)\citenamefont{Schirmer, Fu, and
  Solomon}}]{Schirmer1}
\bibinfo{author}{\bibfnamefont{S.~G.} \bibnamefont{Schirmer}},
  \bibinfo{author}{\bibfnamefont{H.}~\bibnamefont{Fu}}, \bibnamefont{and}
  \bibinfo{author}{\bibfnamefont{A.~I.} \bibnamefont{Solomon}},
  \bibinfo{journal}{Phys. Rev. A} \textbf{\bibinfo{volume}{63}},
  \bibinfo{pages}{063410} (\bibinfo{year}{2001}).

\bibitem[{\citenamefont{Lan et~al.}(2005)\citenamefont{Lan, Tarn, Chi, and
  Clark}}]{Lan1}
\bibinfo{author}{\bibfnamefont{C.}~\bibnamefont{Lan}},
  \bibinfo{author}{\bibfnamefont{T.-J.} \bibnamefont{Tarn}},
  \bibinfo{author}{\bibfnamefont{Q.-S.} \bibnamefont{Chi}}, \bibnamefont{and}
  \bibinfo{author}{\bibfnamefont{J.}~\bibnamefont{Clark}}, \bibinfo{journal}{J.
  Math. Phys.} \textbf{\bibinfo{volume}{46}}, \bibinfo{pages}{052102}
  (\bibinfo{year}{2005}).

\bibitem[{\citenamefont{Wu et~al.}(2006{\natexlab{a}})\citenamefont{Wu, Tarn,
  and Li}}]{wrb1}
\bibinfo{author}{\bibfnamefont{R.-B.} \bibnamefont{Wu}},
  \bibinfo{author}{\bibfnamefont{T.-J.} \bibnamefont{Tarn}}, \bibnamefont{and}
  \bibinfo{author}{\bibfnamefont{C.}~\bibnamefont{Li}}, \bibinfo{journal}{Phys.
  Rev. A} \textbf{\bibinfo{volume}{73}}, \bibinfo{pages}{012719}
  (\bibinfo{year}{2006}{\natexlab{a}}).

\bibitem[{\citenamefont{Puri}(2001)}]{Puri1}
\bibinfo{author}{\bibfnamefont{R.~R.} \bibnamefont{Puri}},
  \emph{\bibinfo{title}{Mathematical methods of quantum optics}}
  (\bibinfo{publisher}{Springer}, \bibinfo{year}{2001}).

\bibitem[{\citenamefont{Gerry and Vrscay}(1989)}]{Gerry2}
\bibinfo{author}{\bibfnamefont{C.~C.} \bibnamefont{Gerry}} \bibnamefont{and}
  \bibinfo{author}{\bibfnamefont{E.~R.} \bibnamefont{Vrscay}},
  \bibinfo{journal}{Phys. Rev. A} \textbf{\bibinfo{volume}{39}},
  \bibinfo{pages}{5717} (\bibinfo{year}{1989}).

\bibitem[{\citenamefont{Gortel and Turski}(1991)}]{Gortel1}
\bibinfo{author}{\bibfnamefont{Z.~W.} \bibnamefont{Gortel}} \bibnamefont{and}
  \bibinfo{author}{\bibfnamefont{L.~A.} \bibnamefont{Turski}},
  \bibinfo{journal}{Phys. Rev. A} \textbf{\bibinfo{volume}{43}},
  \bibinfo{pages}{3221} (\bibinfo{year}{1991}).

\bibitem[{\citenamefont{Bose}(1985)}]{Bose1}
\bibinfo{author}{\bibfnamefont{S.~K.} \bibnamefont{Bose}}, \bibinfo{journal}{J.
  Phys. A} \textbf{\bibinfo{volume}{18}}, \bibinfo{pages}{903}
  (\bibinfo{year}{1985}).

\bibitem[{\citenamefont{Agarwal and Banerji}(2001)}]{Agarwal1}
\bibinfo{author}{\bibfnamefont{G.~S.} \bibnamefont{Agarwal}} \bibnamefont{and}
  \bibinfo{author}{\bibfnamefont{J.}~\bibnamefont{Banerji}},
  \bibinfo{journal}{Phys. Rev. A} \textbf{\bibinfo{volume}{64}},
  \bibinfo{pages}{023815} (\bibinfo{year}{2001}).

\bibitem[{\citenamefont{Wu et~al.}(2006{\natexlab{b}})\citenamefont{Wu, Li, Wu,
  Tarn, and Zhang}}]{wjw1}
\bibinfo{author}{\bibfnamefont{J.-W.} \bibnamefont{Wu}},
  \bibinfo{author}{\bibfnamefont{C.-W.} \bibnamefont{Li}},
  \bibinfo{author}{\bibfnamefont{R.-B.} \bibnamefont{Wu}},
  \bibinfo{author}{\bibfnamefont{T.-J.} \bibnamefont{Tarn}}, \bibnamefont{and}
  \bibinfo{author}{\bibfnamefont{J.}~\bibnamefont{Zhang}}, \bibinfo{journal}{J.
  Phys. A} \textbf{\bibinfo{volume}{39}}, \bibinfo{pages}{13531}
  (\bibinfo{year}{2006}{\natexlab{b}}).

\bibitem[{\citenamefont{Gerry}(1985)}]{Gerry1}
\bibinfo{author}{\bibfnamefont{C.~C.} \bibnamefont{Gerry}},
  \bibinfo{journal}{Phys. Rev. A} \textbf{\bibinfo{volume}{31}},
  \bibinfo{pages}{2721} (\bibinfo{year}{1985}).

\bibitem[{\citenamefont{Walls}(1983)}]{Walls1}
\bibinfo{author}{\bibfnamefont{D.~F.} \bibnamefont{Walls}},
  \bibinfo{journal}{Nature} \textbf{\bibinfo{volume}{306}},
  \bibinfo{pages}{141} (\bibinfo{year}{1983}).

\bibitem[{\citenamefont{Vilenkin and Klimyk}(1991)}]{Vilenkin1}
\bibinfo{author}{\bibfnamefont{N.}~\bibnamefont{Vilenkin}} \bibnamefont{and}
  \bibinfo{author}{\bibfnamefont{A.~U.} \bibnamefont{Klimyk}},
  \emph{\bibinfo{title}{Representation of Lie Groups and Special Functions}},
  vol. \bibinfo{volume}{Volume.1: Simplest Lie groups, special functions and
  integral transforms} (\bibinfo{publisher}{Kluwer Academic Publishers},
  \bibinfo{address}{Boston}, \bibinfo{year}{1991}).

\bibitem[{\citenamefont{Ramakrishna et~al.}(2000)\citenamefont{Ramakrishna,
  Flores, Rabitz, and Ober}}]{Ramakrishna2}
\bibinfo{author}{\bibfnamefont{V.}~\bibnamefont{Ramakrishna}},
  \bibinfo{author}{\bibfnamefont{K.~L.} \bibnamefont{Flores}},
  \bibinfo{author}{\bibfnamefont{H.}~\bibnamefont{Rabitz}}, \bibnamefont{and}
  \bibinfo{author}{\bibfnamefont{R.~J.} \bibnamefont{Ober}},
  \bibinfo{journal}{Phys. Rev. A} \textbf{\bibinfo{volume}{62}},
  \bibinfo{pages}{053409} (\bibinfo{year}{2000}).

\bibitem[{\citenamefont{Perelomov}(1986)}]{Perelomov1}
\bibinfo{author}{\bibfnamefont{P.}~\bibnamefont{Perelomov}},
  \emph{\bibinfo{title}{Generalized Coherent States and Their Applications}}
  (\bibinfo{publisher}{Springer-Verlag}, \bibinfo{address}{Berlin},
  \bibinfo{year}{1986}).

\bibitem[{\citenamefont{Jurdjevic}(1997)}]{Jurdjevic2}
\bibinfo{author}{\bibfnamefont{V.}~\bibnamefont{Jurdjevic}},
  \emph{\bibinfo{title}{Geometric Control Theory}}
  (\bibinfo{publisher}{Cambridge University Press}, \bibinfo{year}{1997}).

\end{thebibliography}
\end{document}